\documentclass{amsart}

\usepackage{amsmath,amssymb,latexsym, amscd}
\usepackage{exscale, cite, eps fig, graphics}

\usepackage{hyperref,xcolor}

\renewcommand{\leq}{\leqslant}

\newtheorem{theorem}{Theorem}[section]

\newtheorem{example}[theorem]{Example}
\newtheorem*{main-theorem}{Main Theorem}
\newtheorem*{remark*}{Remark}
\numberwithin{equation}{section}

\title[Pressure to surface for water waves with vorticity]
{On the recovery of traveling water waves \\ with vorticity from the pressure at the bed} 

\author[Hur]{Vera~Mikyoung~Hur}
\address{Department of Mathematics, University of Illinois at Urbana-Champaign, Urbana, IL 61801 USA}
\email{verahur@math.uiuc.edu}

\author[Livesay]{Michael~R.~Livesay}
\email{mlivesa2@illinois.edu}  

\date{\today}

\begin{document}

\maketitle

\begin{abstract}
We propose higher-order approximation formulae recovering the surface elevation 
from the pressure at the bed and the background shear flow for small-amplitude Stokes and solitary water waves.
They offer improvements over the pressure transfer function and the hydrostatic approximation. 
The formulae compare reasonably well with asymptotic approximations of the exact relation 
between the pressure at the bed and the surface wave in the zero vorticity case,
but they incorporate the effects of vorticity through solutions of the Rayleigh equation. 
Several examples are discussed.
\end{abstract}

\section{Introduction}\label{sec:intro}

A basic problem in oceanography is to determine wave parameters
--- significant wave height, significant wave period, spectral peaks, etc. --- from ocean measurements. 
For instance, the task of tracking the genesis and propagation of tsunamis is of obvious importance.  
One main source of data extensively used for the purpose is 
pressure transducers seeded throughout the Pacific and Indian Ocean.
They collect pressure readings at various water depths and transmit to monitoring stations.

This motivates an interesting mathematical question. Suppose that 
a wave runs in a channel of water over a long distance practically at a constant velocity without change of form, 
and that the value of the pressure at the bed is given, 
and perhaps some other information about the upstream and downstream flow.
From such scant data, can one recover the wave height? 
Incidentally traveling waves may be used as a means to understand more general wave motions; 
see \cite{OVDH2012, deconinck2012relating}, for instance.

Under the assumption that the fluid in the bulk is irrotational, 
a simple approach, which is in practice, for instance, in tsunami detection, is 
to take up the hydrostatic approximation (see \cite{DD1991,KC2010}, for instance):
\begin{equation}\label{def:hydrostatic}
\eta(x)=\frac1gp(x).
\end{equation}
Here $x$ denotes the spatial variable in the direction of wave propagation, 
$\eta$ is the surface displacement from the undisturbed fluid depth $h_0$, say, 
and $p$ is the dynamic pressure, measuring the departure from the hydrostatic pressure;
$g$ is the constant due to gravitational acceleration.
Another is the pressure transfer function (see \cite{DD1991,KC2010}, for instance):
\begin{equation}\label{def:transfer}
\mathcal{F}(\eta)(k)=\frac1g\cosh(kh_0)\mathcal{F}(p)(k).
\end{equation}
Here and throughout,
\[
\mathcal{F}(f)(k)=\int^\infty_{-\infty}f(x)e^{-ikx}~dx
\]
stands for the Fourier transform of the function $x\mapsto f(x)$. 
Note that \eqref{def:transfer} becomes \eqref{def:hydrostatic} in the limit as $h_0\to0$.

Laboratory experiments in \cite{bishop1987measuring}, for instance, support that 
\eqref{def:transfer} satisfactorily predicts the wave height. Furthermore one can derive it consistently
in the regime of small-amplitude Stokes waves in the case of zero vorticity; see \cite{ES2008}, for instance.
On the other hand, numerous studies raised doubts about the adequacy of using the linear theory; 
see \cite{HHK1966, cavaleri1980, biesel1982, LW1984, KC1994}, for instance.
Note that the effects of nonlinearity and current are not negligible in shallow water or in the surf zone; 
see \cite{LW1984}, for instance.

Remarkably, exact relations were derived in \cite{OVDH2012, Constantin2012, CC2013periodic} between 
the trace of the pressure at the horizontal bed and the surface elevation for Stokes and solitary water waves. 
In particular, the formulae apply to large amplitude waves.
They are implicit but, nevertheless, easily implemented in numerical computations, 
and the results agree to varying degrees with laboratory experiments; 
see \cite{deconinck2012relating}, for instance.  

The arguments strongly use that in the case of zero vorticity, 
one is to solve the Cauchy problem for the Laplace equation.
Unfortunately they cannot accommodate underlying shear flows and other physical aspects.
We pause to remark that real flows are hardly irrotational. Rather vorticity is generated, for instance, 
by density stratification, the shear force of the wind, currents or tidal forces, and the effects of bathymetry. 
At present, no exact relations are available between the pressure at the bed and the surface wave in rotational flows. 
Furthermore numerical schemes approximating the exact formulae do not converge, 
because the Cauchy problem for an elliptic PDE is ill-posed.

Recently in \cite{CHW2015}, one of the authors elaborated \eqref{def:transfer} and \eqref{def:hydrostatic}
to permit vorticity and density stratification. 
Specifically, the pressure transfer function and the hydrostatic approximation 
were consistently derived for small-amplitude surface and interface waves 
in an arbitrary shear flow. Unfortunately they do not capture the effects of nonlinearity. 
Furthermore the hydrostatic approximation does not sense the effects of vorticity. 

Here we take matters further and propose higher-order approximation formulae 
recovering the surface elevation from the pressure at the bed 
for small-amplitude Stokes and solitary water waves in an arbitrary shear flow. 
Specifically, we compute higher-order correction terms to 
the pressure transfer function and the hydrostatic approximation in \cite{CHW2015}. 
To the best of the authors' knowledge, these are new. 
We carry out higher-order perturbations of the governing equations, 
rather than relying on a less empirical approach of higher-order Stokes expansion.
We sacrifice the ability to accommodate large amplitude waves, but 
we are able to work to an arbitrary, albeit finite, degree of accuracy, 
when exact formulae relating the pressure at the bed and the surface wave are unavailable. 


The formulae incorporate the effect of vorticity through solutions to the Rayleigh equation, 
which one must in general investigate numerically. But we make an effort to discuss some examples. 
In the case of zero vorticity, in particular, we demonstrate that our results compare reasonably well with 
asymptotic approximations of the exact formulae in \cite{OVDH2012}, for instance;
see Example~\ref{e.g.:periodic0} and Example~\ref{e.g.:solitary0}.
The upshot of the present treatment is highly computationally manageable solutions,
which may develop into an easy and effective numerical scheme.
The practical use of the results, including numerical and experimental studies, is of future investigation.

\section{Preliminaries}\label{sec:preliminaries}
The water wave problem, in the simplest form, concerns the wave motion at the surface 
of an incompressible inviscid fluid, below a body of air and acted upon by gravity. 
For definiteness, we take Cartesian coordinates $(x,y)$, 
where the $x$-axis points in the direction of wave propagation and the $y$-axis vertically upward. 
In other words, the motion is constant in one horizontal direction.
The fluid at time $t$ occupies a region in $\mathbb{R}^2$, bounded 
above by the {\em free} surface and below by the fixed horizontal bottom $y=0$, say. 
Let $y=h(x;t)$ describe the fluid surface at time $t$,
and we assume that $h$ is a single-valued, non-negative and smooth function. 
In the bulk of the fluid, the velocity $(u(x,y;t), v(x,y;t))$ and the pressure $P(x,y;t)$ satisfy 
the Euler equations for an incompressible fluid:
\[\left\{
\begin{split}
&u_t+uu_x+vu_y=-P_x,  \\
&v_t+uv_x+vv_y=-P_y-g \qquad\text{in}\quad 0<y<h(x;t), \\
&u_x+v_y=0.
\end{split}\right.
\]
Here and throughout, a subscript denotes partial differentiation. 
Although an incompressible fluid such as water may have variable density, 
we assume for simplicity that the density $=1$. 
The flow is allowed to be rotational and characterized by the vorticity: 
\[
\omega=v_x-u_y.
\]

The kinematic and dynamic conditions at the fluid surface
\[
v=h_t+uh_x\quad\text{and}\quad P=P_{atm}\qquad \text{at}\quad y=h(x;t)
\]
state, respectively, that each water particle at the surface remains so for all time 
and that the pressure at the surface equals the atmospheric pressure $P_{atm}$; 
we assume that the air is quiescent and we neglect the effects of surface tension.
The flow is required to be tangential to the bottom:
\[
v=0\qquad\text{at}\quad y=0.
\]

It is a matter of experience that waves which are commonly seen in the ocean or a lake propagate 
over a long distance practically at a constant velocity without change of form, namely {\em traveling waves}. 
In other words, $u$, $v$ and $P$ are functions of $(x-ct,y)$ and $h$ is a function of $x-ct$ 
for some $c>0$, the speed of wave propagation. 
Under this assumption, we will go to a moving coordinate frame, changing $x-ct$ to $x$,
whereby the time $t$ completely disappears. The result becomes:
\begin{equation}\label{E:steady}\left\{
\begin{split}
&(u-c)u_x+vu_y=-P_x, && \\
&(u-c)v_x+vv_y=-P_y-g &&\text{in}\quad 0<y<h(x),\\
&u_x+v_y=0, && \\
&v=(u-c)h_x\quad\text{and}\quad P=P_{atm}\qquad &&\text{at}\quad y=h(x),\\
&v=0\qquad&&\text{at}\quad y=0.
\end{split}\right.
\end{equation}

Note that
\begin{equation}\label{def:shear}
h\equiv h_0,\quad(u,v)=(U(y),0)\quad\textrm{and}\quad P =P_{\textrm{atm}}+g(h_0-y)
\end{equation}
make a solution of \eqref{E:steady} for arbitrary $c>0$, $h_0>0$ and an arbitrary $U\in C^1([0,h_0])$. 
Physically, it represents a shear flow, for which 
the velocity and the fluid surface are horizontal and the pressure is hydrostatic. 
The present interest is waves propagating in the $x$-direction over a prescribed shear flow of the form. 
In what follows, therefore, $h_0$ and $U$ are held fixed.
Note that the vorticity of \eqref{def:shear} is $-U'(y)$. 
Here and throughout, the prime denotes ordinary differentiation. 

\subsection*{Scaling of variables}
In order to systematically characterize various approximations, we introduce 
\begin{equation}\label{def:parameters}
\delta=\textrm{the long wavelength parameter}\quad\textrm{and}\quad
\epsilon=\textrm{the amplitude parameter},
\end{equation}
and we define the set of scaled variables.
Rather than introducing a new notation for the variables, we choose, wherever convenient, 
to write, for instance, $x \mapsto x/\delta$. This is to be read that $x$ is replaced by $x/\delta$, 
so that hereafter the symbol $x$ will denote a scaled variable. With this understanding, let 
\begin{equation}\label{def:x-scale}
x\mapsto x/\delta
\end{equation}
and 
\begin{equation}\label{def:uv-scale}
u\mapsto U+\epsilon u_1+\epsilon^2 u_2+\epsilon^3 u_3+\cdots\quad\text{and}\quad 
v\mapsto \delta(\epsilon v_1+\epsilon^2 v_2+\epsilon^3 v_3+\cdots).
\end{equation}
Moreover, we write
\begin{equation}\label{def:h-scale}
h=h_0+\epsilon\eta_1+\epsilon^2\eta_2+\epsilon^3\eta_3+\cdots
\end{equation}
and 
\begin{equation}\label{def:P-scale}
P=P_{atm}+g(h_0-y)+\epsilon p_1+\epsilon^2 p_2+\epsilon^3 p_3+\cdots.
\end{equation}
Physically, $h-h_0$ means the surface displacement from the undisturbed fluid depth 
and $P-P_{atm}-g(h_0-y)$ is the dynamic pressure, measuring the departure from the hydrostatic pressure;
see \cite[Section~1.3.2 and Section~3.4.1]{Johnson1997} for the detail.

Substituting \eqref{def:x-scale} through \eqref{def:P-scale} into \eqref{E:steady} 
and restricting the result to the undisturbed fluid domain $0<y<h_0$, we arrive at that:
\begin{equation}\label{E:Euler-scaled}
\left\{\begin{split}
&(U-c)u_{1x}+U'v_1+\epsilon((U-c)u_{2x}+U'v_2+u_1u_{1x}+v_1u_{1y})\\ 
&+\epsilon^2((U-c)u_{3x}+U'v_3+(u_1u_{2})_x+u_2u_{1x}+v_1u_{2y}+v_2u_{1y})+\cdots\\
&\hspace*{185pt}=-p_{1x}-\epsilon p_{2x}-\epsilon^2 p_{3x}+\cdots,	\\
&\delta^2((U-c)v_{1x}+\epsilon((U-c)v_{2x}+u_1v_{1x}+v_1v_{1y})\\
&+\epsilon^2((U-c)v_{3x}+u_1v_{2x}+u_2v_{1x}+v_1v_{2y}+v_2v_{1y}))+\cdots\\
&\hspace*{185pt}=-p_{1y}-\epsilon p_{2y}-\epsilon^2 p_{3y}+\cdots,\\
&u_{1x}+v_{1y}+\epsilon(u_{2x}+v_{2y})+\epsilon^2(u_{3x}+v_{3y})+\cdots=0
\end{split}\right.
\end{equation}
in $0<y<h_0$, and
\begin{equation}\label{E:top-scaled}
\left\{\begin{split}
&(U-c)\eta_{1x}+\epsilon((U-c)\eta_{2x}+u_1\eta_{1x})\\
&+\epsilon^2((U-c)\eta_{3x}+u_1\eta_{2x}+u_2\eta_{1x})+\cdots
=v_1+\epsilon v_2+\epsilon^2 v_3+\cdots,\\
&g(\eta_1+\epsilon \eta_2+\epsilon^2\eta_3)+\cdots
=p_1+\epsilon p_2+\epsilon^2 p_3+\cdots  
\end{split}\right.
\end{equation}
at $y=h_0$ and
\begin{equation}\label{E:bottom-scaled}
v_1+\epsilon v_2+\epsilon^2 v_3+\cdots= 0 \qquad\text{at}\quad y=0.
\end{equation}

Note that $u=v=p=\eta=0$ --- no disturbances --- satisfy \eqref{E:Euler-scaled}-\eqref{E:bottom-scaled}
for arbitrary $c>0$, $h_0>0$ and an arbitrary $U\in C^1([0,h_0])$.

The governing equations may be non-dimensionalized (see \cite[Section~1.3.1]{Johnson1997}, for instance),
but we do not pursue it here. Rather we work in the physical variables as much as possible.

\section{Stokes waves}\label{sec:Stokes}
By {\em Stokes waves}, we mean solutions of \eqref{E:steady}, which are periodic and symmetric in the $x$-direction
(historically, in the case of zero vorticity and practically at rest at great depths).
For an arbitrary distribution of vorticity, under some assumptions, incidentally, 
all periodic solutions of \eqref{E:steady} are a priori symmetric about their crests; 
see \cite{Hur2007symm} and \cite{CEW2007symm}, for instance. 
Stokes in his classic memoir in 1847 (see also \cite{Stokes}) made many contributions about waves of the kind, 
observing, for instance, that crests tend to be sharper and troughs flatter as the amplitude increases 
and that the crest of a wave of greatest height would be a stagnation point with a $120^\circ$ corner. 

In the case of zero vorticity, the rigorous existence theory of Stokes waves goes back to constructions
in \cite{Nekrasov, Levi-Civita} and \cite{Struik} of small amplitude waves,
and it includes global bifurcation results in \cite{Krasovskii, KN}, for instance,
and the resolution in \cite{AFT} of Stokes' conjecture about the wave of greatest height. 
All these works strongly use the assumption that the flow in the bulk is irrotational,
whereby one may reformulate the problem in terms of quantities defined at the fluid surface. 
We encourage the interested reader to some excellent surveys 
\cite{Toland1996Stokes, OS2001book, BTbook, Groves2004survey, Strauss2010}. 

The zero vorticity assumption may serve as a reasonable approximation in some circumstances.
Moreover in the absence of initial vorticity, boundaries or currents, 
water waves will have zero vorticity for all future time.
But rotational effects are significant in many circumstances, for instance,
for wind-driven waves, waves riding upon a sheared current, or waves near a ship or pier. 

In the case of nonzero vorticity, the situation is to look inside the fluid 
because the velocity potential is no longer viable to use.
Consequently it is harder to handle, analytically and numerically, than the zero vorticity case.
It was not until recently that Constantin and Strauss \cite{CS2004} established the existence, 
from zero up to (but not including) an ``extreme" wave exhibiting a stagnation point. 
Specifically, for arbitrary wave speed and period and 
for an arbitrary function relating the vorticity and the stream function, 
subject to a ``bifurcation condition", they constructed a global continuum of Stokes waves with vorticity. 
This quickly led to a flurry of research activities. 
It would be impossible to do justice to all the advances in the direction, but we single out a few ---
\cite{Hur2006Stokes1, Hur2011Stokes2} in the infinite depth, 
\cite{Wahlen2009critical, EEW2011critical, CSV2014} permitting critical layers, 
\cite{Varvaruca2009extreme, VW2012extreme} about an extreme wave, 
and \cite{constantin2011discontinuous} permitting discontinuous vorticities.

A key idea in \cite{CS2004}, as in all free boundary problems, is to fix the --- a priori unknown --- fluid domain. 
As a matter of fact, the bifurcation condition in \cite{CS2004} asks 
if the linearization of \eqref{E:steady} about \eqref{def:shear} admits a nontrivial solution, 
but in the Dubreil-Jacotin variables, which map the fluid domain of one period to a fixed rectangle. 
Seeking explicit relations between the pressure at the bed and the surface wave in the physical coordinates 
(where the underlying shear flow and the fluid depth are fixed), 
here we carry out all calculations in the physical variables. 
Below we record the translation of the bifurcation condition in \cite{CS2004} to the physical variables, 
originally derived by one of the authors in \cite{HL2008}.  

\subsection*{The bifurcation condition}
For arbitrary $h_0>0$ and $U\in C^1([0,h_0])$, 
recall that \eqref{def:shear} is a solution of \eqref{E:steady} for all $c>0$. 
We are interested in determining at which values of $c>0$ and $k>0$, 
there bifurcates a family of small-amplitude $2\pi/k$-periodic solutions of \eqref{E:steady}. 
A necessary condition, it turns out, is that 
\begin{equation}\label{E:bifur}\left\{
\begin{split}
&(U-c)(\phi''-k^2\phi)-U''\phi=0\qquad\textrm{for }0<y<h_0,\\
&\phi'(h_0)=\Big(\frac{g}{(U(h_0)-c)^2}+\frac{U'(h_0)}{U(h_0)-c}\Big)\phi(h_0)
\quad\textrm{and}\quad\phi(0)=0
\end{split}\right.
\end{equation}
admits a nontrivial solution for some $c$ and $k$.
It is in general not a sufficient condition, but in case when ${\displaystyle c>\max_{0\leq y\leq h_0}U(y)}$,
bifurcation does occur, provided that the kernel of \eqref{E:bifur} is one dimensional.
Under this assumption, furthermore, $u<c$ throughout the fluid region. Note that $u=c$ at a stagnation point.
In \cite{CS2004}, the wave speed and the vorticity-stream function relation are fixed 
whereas the shear flow $U$ and the fluid depth $h_0$ vary along the branch of solutions. 
On the contrary, here the shear flow and the fluid depth are fixed 
and the wave speed is determined upon solving \eqref{E:bifur}. 
The ordinary differential equation in \eqref{E:bifur} 
goes by the name of the Rayleigh (or inviscid Orr-Sommerfeld) equation. 
It is not singular, provided that $c > \max U$. 

In the case of $U\equiv0$, namely the zero vorticity, 
a straightforward calculation reveals that a nontrivial solution to \eqref{E:bifur} exists, provided that
\begin{equation}\label{E:c0}
c^2=\frac{g\tanh(kh_0)}{k}.
\end{equation}
This is the well-known dispersion relation of water waves in irrotational flows.  
In the case of $U(y)=\gamma y$ for some constant $\gamma$, namely the constant vorticity $-\gamma$, 
similarly, a straightforward calculation reveals that a nontrivial solution to \eqref{E:bifur} exists, provided that  
\begin{equation}\label{E:cgamma}
c=-\frac{\gamma\tanh(kh_0)}{2k}+\sqrt{\frac{\gamma^2\tanh^2(kh_0)}{4k^2}+\frac{gk\tanh(kh_0)}{k}}.
\end{equation}
This is the dispersion relation in the case of constant vorticity.
(The other solution with the $-$ sign violates $c>\max U$, and hence we discard it.)

The bifurcation condition is closely related, but not equivalent, to the dispersion relation. 
The bifurcation condition is a necessary and sufficient condition for nontrivial solutions to exist
whereas the dispersion relation is a necessary condition 
for plane wave solutions to exist to the associated linear problem.

For general non-constant vorticities, one must not expect to solve \eqref{E:bifur} explicitly; 
see \cite{Pete2012dispersion}, for instance.
For a wide range of shear flows, nevertheless, one may be able to verify the bifurcation condition using the ODE theory. 
If $U\in C^2([0,h_0])$, $U''(h_0)<0$ and $U(h_0)>U(y)$ for $0<y<h_0$, for instance, 
bifurcation takes place for some $c>\max U$ for all $k>0$; see \cite[Lemma~2.5]{HL2008}. 
In the long wave limit as $k\to0+$, suitable for solitary water waves, 
the bifurcation condition leads to the Burns condition
\begin{equation}\label{E:Burns}
\int^{h_0}_0 \frac{dy}{(U(y)-c)^2}=\frac1g;
\end{equation}
see \cite{Burns, Johnson1997, HL2008}, for instance.

\subsection{The first-order approximation}\label{sec:Stokes1}
We set forth the surface reconstruction procedure from the pressure at the fluid bed 
for small-amplitude Stokes waves with vorticity. We therefore assume that 
\begin{equation}\label{E:Stokes-parameter}
\delta=1\quad\text{and}\quad\epsilon\ll1
\end{equation}
and, without loss of generality, 
$u_j$'s, $v_j$'s, $p_j$'s and $\eta_j$'s, $j=1,2,3, \dots$, in \eqref{def:uv-scale}-\eqref{def:P-scale}
are $2\pi$-periodic in the $x$-variable. In other words, the wave number $k=1$. 
We do not assume a priori their symmetry and monotonicity.

Under this assumption, 
\eqref{E:Euler-scaled} and \eqref{E:top-scaled} \eqref{E:bottom-scaled} at the leading order become:
\begin{equation}\label{E:periodic-Euler1}\left\{
\begin{split}
&(U-c)u_{1x}+U'v_1= -p_{1x},\quad&&\\
&(U-c)v_{1x}= -p_{1y}\qquad&&\text{in}\quad 0<y<h_0, \\
&u_{1x}+v_{1y}= 0,&&\\
\end{split}
\right.\end{equation}
and
\begin{align}
&v_1=(U-c)\eta_{1x}\quad\text{and}\quad p_1=g\eta_1&&\text{at}\quad y=h_0, \label{E:periodic-top1}\\
&v_1=0 &&\text{at}\quad y=0.\label{E:periodic-bottom1}
\end{align}

Stokes waves with vorticity constructed in \cite{CS2004}, for instance, solve \eqref{E:steady}, 
and hence \eqref{E:Euler-scaled}-\eqref{E:bottom-scaled}.
It follows from bifurcation theory that small-amplitude solutions solve 
\eqref{E:periodic-Euler1}-\eqref{E:periodic-bottom1} with errors in $O(\epsilon^2)$ in some suitable H\"older space. 
Thanks to the symmetry of the bifurcation problem, furthermore, 
the wave speed is approximated by that determined upon solving \eqref{E:bifur} with errors in $O(\epsilon^3)$.
Throughout the section $c>0$  is held fixed.

Differentiating the first equation in \eqref{E:periodic-Euler1} in the $y$-variable 
and the second equation in the $x$-variable, and using the third equation, we arrive at that 
\[
(U-c)\Delta v_1-U''v_1=0 \qquad \text{in}\quad 0<y<h_0.
\]
Suppose that the dynamic pressure (see \eqref{def:P-scale}) at the bed is prescribed to the order of $\epsilon$. 
We assume that it is smooth and $2\pi$-periodic in the $x$-variable, 
and thus we write it in the (complex) Fourier series as
\begin{equation}\label{def:periodic-b1s}
p_1(x,0)=\sum_{n=-\infty}^\infty b_{1n}e^{inx}.
\end{equation}
The first and third equations in \eqref{E:periodic-Euler1} restricted at the bed imply that
\[
(U(0)-c)v_{1y}(x,0)-U'(0)v_1(x,0)=p_{1x}(x,0)=\sum_{n=-\infty}^\infty inb_{1n}e^{inx}.
\]
Note from \eqref{E:periodic-bottom1} that the second term on the left side vanishes. 
To recapitulate, in the first-order approximation of the surface elevation as a function of the pressure at the bed
for small-amplitude Stokes waves with vorticity, one solves the Cauchy problem for the linear elliptic PDE:
\begin{equation}\label{E:periodic-prob1}
\left\{\begin{split}
&(U-c)\Delta v_1-U''v_1=0 \qquad &&\text{in}\quad 0<y<h_0,\\
&v_1=0 &&\text{at}\quad y=0,\\
&v_{1y}=\frac{1}{U(0)-c}\sum_{n=-\infty}^\infty inb_{1n}e^{inx}\qquad&&\text{at}\quad y=0,
\end{split}\right.
\end{equation}
and one uses equations in \eqref{E:periodic-Euler1} 
to determine $u_1$ and $p_1$ in terms of $p_1(\cdot,0)$ (see \eqref{def:periodic-b1s}). 
One then determines $\eta_1$ using the second equation in \eqref{E:periodic-top1} as $g\eta_1(x)=p_1(x,h_0)$.

Specifically, let's write that
\[
v_1(x,y)=\sum_{n=-\infty}^\infty i\phi_{1n}(y)e^{inx}.
\]
Here $i$  is for convenience.
For each $n$, note that \eqref{E:periodic-prob1} leads to the Cauchy problem for the Rayleigh equation:
\begin{equation}\label{E:phi1n}
\left\{\begin{split}
&(U-c)(\phi_{1n}''-n^2\phi_{1n})-U''\phi_{1n}=0\qquad \text{for}\quad 0<y<h_0,\\
&\phi_{1n}(0)= 0,\\
&\phi_{1n}'(0)= \frac{nb_{n}}{U(0)-c}.
\end{split}\right.
\end{equation}
It follows from the ODE theory that \eqref{E:phi1n} admits a unique solution. In particular $\phi_{10}\equiv0$. 
We then infer from the last equation in \eqref{E:periodic-Euler1} that 
\[
u_1(x,y)= u_{10}(y)-\sum_{n\neq 0}\phi_{1n}'(y)\frac{e^{inx}}{n}
\]
for some function $u_{10}$ and, similarly, from the first equation in \eqref{E:periodic-Euler1} that
\[
p_1(x,y)= p_{10}(y)+\sum_{n\neq 0}((U(y)-c)\phi_{1n}'(y)-U'(y)\phi_{1n}(y))\frac{e^{inx}}{n}
\]
for some function $p_{10}$. 
Comparing $p_{1y}$ to the second equation in \eqref{E:periodic-Euler1}, 
we use \eqref{E:phi1n} to deduce that $p_{10}(y)$ is a constant. 
Comparing this to \eqref{def:periodic-b1s}, moreover, we conclude that $p_{10}(y)=b_0$. 
One may not be able to determine $u_{10}$, on the other hand. This is not surprising 
since adding a function of $y$ to $u_1$ does not change \eqref{E:periodic-Euler1}-\eqref{E:periodic-bottom1}. 
Ultimately the second equation in \eqref{E:periodic-top1} implies that 
\begin{equation}\label{E:periodic1s}
g\eta_1(x)=b_0+\sum_{n\neq 0}((U(h_0)-c)\phi_{1n}'(h_0)-U'(h_0)\phi_{1n}(h_0))\frac{e^{inx}}{n}.
\end{equation}
This furnishes an implicit formula relating $\eta_1$ to $p_1(\cdot,0)$ (see \eqref{def:periodic-b1s}),
provided with the background shear flow and the wave speed, subject to the bifurcation condition. 
It incorporates the effects of vorticity through the solution to the Cauchy problem for the Rayleigh equation.

In what follows, we furthermore assume that 
\begin{equation}\label{assume:Stokes}
\text{$\eta_1$ and $u_1(\cdot, y)$ for each $y\in [0,h_0]$ are proportional to $\cos(x)$.}
\end{equation}
It follows from local bifurcation theory that the leading part of small-amplitude solutions 
constructed in \cite{CS2004}, for instance, satisfy \eqref{assume:Stokes}. 
As a matter of fact, all solutions in \cite{CS2004} satisfy \eqref{assume:Stokes}, 
after possibly redefining the undisturbed fluid depth $h_0$ and the background current $U$.
Under this assumption, the pressure at the bed must be prepared as 
\begin{equation}\label{def:periodic-b1}
p_1(x,0)=b\cos(x)
\end{equation}
for some constant $b$ ($b=2b_{11}=2b_{1-1}$), and a straightforward calculation reveals that 
\begin{equation}\label{E:periodic-soln1}
\left\{\begin{split}
&u_1(x,y)=\phi'(y)\cos(x),\\
&v_1(x,y)=\phi(y)\sin(x), \\
&p_1(x,y)=((c-U(y))\phi'(y)+U'(y)\phi(y))\cos(x)
\end{split}\right.
\end{equation}
and
\begin{equation}\label{E:periodic-eta1}
\eta_1(x)=\frac1g((c-U(h_0))\phi'(h_0)+U'(h_0)\phi(h_0))\cos(x),
\end{equation}
where $\phi$ is the unique solution of 
\begin{equation}\label{E:Rayleigh1}
\left\{\begin{split}
&(U-c)(\phi''-\phi)-U''\phi=0\qquad \text{for}\quad 0<y<h_0,\\
&\phi(0)= 0,\\
&\phi'(0)= \frac{b}{c-U(0)}.
\end{split}\right.
\end{equation}
This furnishes an implicit formula relating $\eta_1$ and $p_1(\cdot,0)$ (see \eqref{def:periodic-b1}), 
assuming that small-amplitude Stokes waves are sinusoidal to the leading order. 
Note that \eqref{assume:Stokes} uniquely determines $u_1$. 

\begin{remark*}[The pressure transfer function]\rm
The last equation in \eqref{E:periodic-soln1} is reminiscent of the pressure transfer function found in \cite{CHW2015}:
\begin{equation}\label{E:transfer}
p(x,y)=((c-U(y))\phi'(y)+U'(y)\phi(y))\eta(x),
\end{equation}
but $\phi$ in \cite{CHW2015} solves the boundary value problem for the Rayleigh equation: 
\[
\left\{\begin{split}
&(U-c)(\phi''-\phi)-U''\phi=0\qquad \text{for}\quad 0<y<h_0,\\
&\phi(0)= 0,\\
&\phi(h_0)=c-U(h_0).
\end{split}\right.
\]
One may use \eqref{E:transfer} to approximate the pressure at the bed as a function of the surface elevation
for small-amplitude Stokes waves with vorticity.
\end{remark*}

\subsection{The second-order approximation}\label{sec:Stokes2}
We continue to assume \eqref{E:Stokes-parameter}. To proceed,
\eqref{E:Euler-scaled} and \eqref{E:top-scaled}, \eqref{E:bottom-scaled} at the order of $\epsilon$ become:
\begin{equation}\label{E:periodic-Euler2}\left\{
\begin{split}
&(U-c)u_{2x}+U'v_2+u_1u_{1x}+v_1u_{1y}= -p_{2x}, \\
&(U-c)v_{2x}+u_1v_{1x}+v_1v_{1y}= -p_{2y}	&\quad\text{in}\quad 0<y<h_0, \\
&u_{2x}+v_{2y}	= 0, \\
\end{split}
\right.\end{equation}
and
\begin{align}
&v_2=(U-c)\eta_{2x}+u_1\eta_{1x}\quad\text{and}\quad p_2=g\eta_2&&\text{at}\quad y=h_0, \label{E:periodic-top2}\\
&v_2=0 &&\text{at}\quad y=0.\label{E:periodic-bottom2}
\end{align}
Furthermore we assume \eqref{assume:Stokes}, 
and thus $u_1$, $v_1$ and $p_1$ are computed using \eqref{E:periodic-soln1};
otherwise, expressions become quite complicated involving double Fourier series.

Suppose that the dynamic pressure at the bed is prescribed to the order of $\epsilon^2$.
We continue to assume that it is smooth and $2\pi$-periodic in the $x$-variable, and we write that  
\begin{equation}\label{def:periodic-b2s}
p_2(x,0)=\sum_{n=-\infty}^{\infty}b_{2n}e^{inx}.
\end{equation}
We then repeat the argument in the previous subsection.
To summarize, in the second-order approximation of the surface elevation as a function of the pressure at the bed 
for small-amplitude Stokes waves with vorticity,
one solves the Cauchy problem for the inhomogeneous linear elliptic PDE:
\begin{equation}\label{E:periodic-prob2}\left\{
\begin{split}
&(U-c)\Delta v_2-U''v_2=\frac{1}{2}(\phi\phi'''-\phi'\phi'')\sin(2x)\qquad&&\text{in}\quad 	0<y<h_0,\\
&v_2= 0 &&\text{at}\quad y=0,\\
&v_{2y}= \frac{1}{U(0)-c} \Big(\frac{1}{2}(\phi\phi''-(\phi')^2)\sin(2x)
+\sum\limits_{-\infty}^{\infty}inb_{2n}e^{inx}\Big)&&\text{at}\quad y=0,
\end{split}
\right. \end{equation}
and uses equations in \eqref{E:periodic-Euler2} and \eqref{E:periodic-soln1}  to determine $u_2$ and $p_2$ 
in terms of $p_1(\cdot,0)$ and $p_2(\cdot,0)$ (see \eqref{def:periodic-b1} and \eqref{def:periodic-b2s}).
One then determines $\eta_2$ using the second equation in \eqref{E:periodic-top2} as $g\eta_2(x)=p_2(x,h_0)$.

Let's write that 
\[
v_2(x,y)=\sum\limits_{n=-\infty}^{\infty}i\phi_{2n}(y)e^{inx}.
\]
Here $i$ is for convenience. 
For each $n\neq \pm 2$, note that \eqref{E:periodic-prob2} leads to the Cauchy problem for the Rayleigh equation:
\begin{equation}\label{E:phi2n}\left\{
\begin{split}
&(U-c)(\phi_{2n}''-n^2\phi_{2n})-U'' \phi_{2n}= 0 \qquad \text{for}\quad 0<y<h_0,\\
&\phi_{2n}(0)= 0,\\
&\phi_{2n}'(0)= \frac{nb_{2n}}{U(0)-c}.
\end{split}
\right. \end{equation}
For $n=\pm2$, similarly, 
\begin{equation}\label{E:phi22}\left\{
\begin{split}
&(U-c)(\phi_{2\pm2}''-4\phi_{2\pm2})-U'' \phi_{2\pm2}=\mp\frac{1}{4}(\phi\phi'''-\phi'\phi'')e^{\pm2ix}
\quad \text{for }0<y<h_0,\\
&\phi_{2\pm2}(0)= 0, \\
&\phi_{2\pm2}'(0)=\frac{1}{U(0)-c}\Big(\mp\frac{1}{4}(\phi\phi''-(\phi')^2)(0)\pm2b_{2\pm2}\Big).
\end{split}
\right. \end{equation}
It follows from the ODE theory that \eqref{E:phi2n} and \eqref{E:phi22} admit unique solutions. 
In particular $\phi_{20}\equiv0$. We then infer from the last equation in \eqref{E:periodic-Euler2} that 
\[
u_2(x,y)= u_{20}(y)-\sum_{n\neq 0}\phi_{2n}'(y)\frac{e^{inx}}{n}
\]
for some function $u_{20}$ and, similarly, from the first equation in \eqref{E:periodic-Euler2} that
\[
p_2(x,y)=\frac{1}{4}(\phi\phi''-(\phi')^2)(y)\cos(2x)
+p_{20}(y)+\sum_{n\neq 0}((U(y)-c)\phi_{2n}'(y)-U'(y)\phi_{2n}(y))\frac{e^{inx}}{n}
\]
for some function $p_{20}$. Comparing $p_{2y}$ to the second equation in \eqref{E:periodic-Euler2}, 
we use \eqref{E:phi2n} to deduce that $p_{20}'(y)=-(\phi\phi')(y)$.  
Comparing this to \eqref{def:periodic-b2s}, moreover, we conclude that $p_{20}(y)=-\frac{1}{2}\phi^2(y)+b_{20}$.  
One may not able to determine $u_{20}$, on the other hand. This is not surprising 
since adding a function of $y$ to $u_2$ does not change \eqref{E:periodic-Euler2}-\eqref{E:periodic-bottom2}.
Ultimately the second equation in \eqref{E:periodic-top2} implies that 
\begin{equation}\label{E:periodic-eta2}
\begin{split}
g\eta_2(x)=&-\frac{1}{2}\phi^2(h_0)+\frac{1}{4}(\phi\phi''-(\phi')^2)(h_0)\cos(2x)\\
&+b_{20}+\sum\limits_{n\neq 0}((U(h_0)-c)\phi_{2n}'(h_0)-U'(h_0)\phi_{2n}(h_0))\frac{e^{inx}}{n},
\end{split}
\end{equation}
where $\phi$ solves \eqref{E:Rayleigh1} and $\phi_{2n}$ solves \eqref{E:phi2n} or \eqref{E:phi22}.
This furnishes an implicit formula relating $\eta_2$ to $p_1(\cdot,0)$ and $p_2(\cdot,0)$ 
(see \eqref{def:periodic-b1} and \eqref{def:periodic-b2s}), 
provided with the background shear flow and the wave speed, subject to the bifurcation condition. 
It incorporates the effects of vorticity through solutions of the Cauchy problems for the Rayleigh equations.
Note that \eqref{E:periodic-eta2} depends linearly on the second-order pressure data and nonlinearly on the first-order pressure data.

We may repeat the above argument and continue to higher-order approximations.
Expressions become quite complicated, however, and hence we do not pursue here. 
One shortcoming of our method is that we only determine the horizontal velocities $u_j$'s, $j=1,2,3,\dots$, 
up to functions of $y$. We may continue to assume that (see \eqref{assume:Stokes})
\[
\int^{2\pi}_0u_j(x,y)~dx=0\qquad\text{for each}\quad y\in[0,h_0],
\]
$j=1,2,3,\dots$, after possibly redefining the shear flow, to uniquely determine them. 
In the case of zero vorticity, in \cite{OVDH2012}, for instance, 
likewise, one determines the velocity potential up to a constant, 
but the reconstruction formula does not require knowledge of the velocity potential itself; 
see Appendix~\ref{sec:OVDH2012} for the detail.

Another shortcoming is that the wave speed agrees with the bifurcation speed up to the order of $\epsilon^2$
in the regime of small amplitude waves, even in the case of zero vorticity. 
Moreover it is in practice difficult to measure.
We may continue to assume the bifurcation speed for higher-order approximations. 
As a matter of fact, numerical computations in \cite{deconinck2012relating}, for instance, indicate that 
the maximum wave height does not suffer much from the small amplitude limit of the wave speed.

\subsection{Examples}\label{sec:periodic.e.g.}

Formulae in Section~\ref{sec:Stokes1} and Section~\ref{sec:Stokes2} 
must in general be investigated numerically.
In some cases, nevertheless, analytical solutions are available, as we discuss below. 
Throughout the subsection, we assume for simplicity that 
\begin{equation}\label{assume:pStokes}
p_1(x,0)=b\cos(x)\quad\text{and}\quad p_2(x,0)=0.
\end{equation}

\begin{example}[Zero vorticity]\label{e.g.:periodic0}\rm
In the case of $U\equiv0$, namely the zero vorticity, a straightforward calculation reveals that 
the unique solution of \eqref{E:Rayleigh1} is
\[
\phi(y)=\frac{b}{c}\sinh(y),
\]
whence \eqref{E:periodic-soln1} and \eqref{E:periodic-eta1} become
\[
\left\{\begin{split}
&u_1(x,y)=\frac{b}{c}\cosh(y)\cos(x),\\
&v_1(x,y)=\frac{b}{c}\sinh(y)\sin(x),\\
&p_1(x,y)=b\cosh(y)\cos(x),
\end{split}\right.
\]
and 
\begin{equation}\label{E:periodic0eta1}
g\eta_1(x)=b\cosh(h_0)\cos(x).
\end{equation}
This agrees with the pressure transfer function in \eqref{def:transfer} in the case when the wave number $k=1$. 

To proceed, note that $\phi\phi''-(\phi')^2=-(b/c)^2$.  A straightforward calculation reveals that 
the unique solution of \eqref{E:phi2n} or \eqref{E:phi22} is
\[
\phi_{2n}(y)\equiv0 \quad\text{if }n\neq \pm 2
\quad\text{and}\quad \phi_{22}(y)=-\phi_{2-2}(y)=-\frac{b^2}{8c^3}\sinh(2y).
\]
Therefore \eqref{E:periodic-eta2} becomes 
\begin{equation}\label{E:periodic0eta2}
g\eta_2(x)=\frac{1}{2}\left(\frac{b}{c}\right)^2\sinh^2(h_0)(-1+\cos(2x)).
\end{equation}

We compare \eqref{E:periodic0eta1} and \eqref{E:periodic0eta2} to asymptotic approximations of the exact formula 
in \cite{OVDH2012}, for instance, but in the physical variables. 
We include the detail of the algebra in Appendix~\ref{sec:OVDH2012} for completeness. 
The first-order approximation of the result in \cite{OVDH2012} agrees with \eqref{E:periodic0eta1};
see \eqref{E:OVDH1'}.
The second-order approximation of the result in \cite{OVDH2012} may be written, abusing notation, as
\begin{equation}\label{E:eta2OVDH}
g\eta_{2}(x)= \left(\frac{b}{2c} \right)^2 (1+4\sinh^2(h_0) \cos(2x)).
\end{equation}
This shares with \eqref{E:periodic0eta2} that the second-harmonic correction becomes negligible as $h_0\to 0$.
But in \eqref{E:periodic0eta2} the second-order depth correction decreases as $h_0\to0$
whereas in \eqref{E:eta2OVDH} it increases with small depths. 
Note from \eqref{E:c0} that $c\to \sqrt{gh_0}$ as $h_0\to 0$.

A potential reason for disagreements between \eqref{E:periodic0eta2} and \eqref{E:eta2OVDH} 
is that equations in \eqref{E:top-scaled}, and hence \eqref{E:periodic-top2}, 
do not take into account of the nonlinear effects of the boundary conditions at the free surface. 
If we were to include boundary variations in the derivation of \eqref{E:top-scaled} 
(see \cite{HL2008}, for instance) then perhaps we would be able to find a better agreement.

Furthermore we compare the results to the small amplitude asymptotics 
of a true solution in \cite[Section~13.13]{Whitham}, for instance.
The first-order approximation agrees with \eqref{E:periodic0eta1} and \eqref{E:OVDH1'}. 
The second-order approximation of the result in \cite{Whitham} may be written, abusing notation, as
\begin{equation}\label{E:eta2exact}
g\eta_{2}(x)=\frac14\Big(\frac{b}{c}\Big)^2\Big(-\frac{1}{\cosh^2(h_0)}+\Big(2+\frac{3}{\sinh^2(h_0)}\Big)\cos(2x)\Big).
\end{equation}
Observe that the nonlinear corrections become significant as $h_0\to0$.
Experimental studies in \cite{LW1984}, for instance, bear out this.
(Note that the result in \cite{Whitham} assumes that 
the integral of $\eta_j$ over one period is zero for all $j=1,2,3,\dots$.
One may think of this as a consequence of absorbing a constant of integration in the hydrostatic pressure. 
It is straightforward to keep track of nonzero mean values in $\eta_j$'s, however. We omit the detail.) 

A potential reason for disagreements between \eqref{E:periodic0eta2}, \eqref{E:eta2OVDH} and \eqref{E:eta2exact} 
is that \eqref{assume:pStokes} may not hold for the pressure distribution of a true solution. 
If we were to prescribe a more physically realistic pressure input, 
then perhaps we would be able to find a better agreement. 
It is interesting to find a higher-order approximation formula, 
for which the nonlinear effects become significant for small depths.
\end{example}

\begin{example}[Constant vorticity]\rm
Let $U(y)=\gamma y$, $0\leq y\leq h_0$, for some constant $\gamma$. 
(More generally, one may take $U(y)=\gamma y+U_0$ for some constant $U_0$, 
but $U_0$ may be absorbed into the wave speed.) 
This models the constant vorticity $-\gamma$. 

A straightforward calculation reveals that the unique solution of \eqref{E:Rayleigh1} is 
\[
\phi(y)=\frac{b}{c}\sinh(y),
\]
whence \eqref{E:periodic-soln1} and \eqref{E:periodic-eta1} become
\[
\left\{\begin{split}
&u_1(x,y)=\frac{b}{c}\cosh(y)\cos(x),\\
&v_1(x,y)=\frac{b}{c}\sinh(y)\sin(x),\\
&p_1(x,y)=\Big(\frac{c-\gamma y}{c}\cosh(y)+\frac{\gamma}{c}\sinh(y)\Big)b\cos(x)
\end{split}\right.
\]
and 
\[
g\eta_1(x)=\left(\frac{c-\gamma h_0}{c}\cosh(h_0)+\frac{\gamma}{c}\sinh(h_0)\right)b\cos(x).
\]

To proceed, note that $\phi\phi''-(\phi')^2=-(b/c)^2$. A straightforward calculation reveals that 
the unique solution of \eqref{E:phi2n} or \eqref{E:phi22} is
\[
\phi_{2n}(y)\equiv0\quad\text{if }n\neq\pm 2 \quad\text{and}\quad
\phi_{22}(y)=-\phi_{2-2}(y)=-\frac{b^2}{8c^3}\sinh(2y).
\]
Therefore \eqref{E:periodic-eta2} becomes
\begin{align*}
g\eta_2(x)=\frac14\Big(\frac{b}{c}\Big)^2\Big(-2\sinh^2(h_0)
+\Big(\frac{c-\gamma h_0}{c}\cosh(2h_0)+\frac{\gamma}{2c}\sinh(2h_0)-1\Big)\cos(2x)\Big).
\end{align*}
In the case of $\gamma=0$, these formulae reduce to those in Example~\ref{e.g.:periodic0}.
\end{example}

\section{Solitary water waves}\label{sec:solitary}
By {\em solitary water waves}, we mean solutions of \eqref{E:steady}, 
for which instead of the periodic boundary condition, 
$h(x)$ tends to a constant and $v(x,y)\to 0$ uniformly for $y$ as $x\to\pm\infty$. 
Historically they have stimulated a considerable part of developments in the theory of wave motion,
from Russell's famous horseback observations to the elucidation of the Korteweg-de Vries (KdV) solitons.

Solitary water waves may formally be viewed as the limit of Stokes waves as the period tends to infinity. 
As a matter of fact, small-amplitude solitary water waves, if exist, emanate  
near the critical speed determined upon solving \eqref{E:Burns}. 
They are a genuinely nonlinear phenomenon, however, and classical bifurcation theory fail to yield the existence. 
In the case of zero vorticity, their rigorous theory goes back to 
constructions in \cite{FH1954, Beale1977solitary} of small amplitude waves 
and it includes a large amplitude result in \cite{AT1981solitary}. 
Recently, these results have been extended in the case of non-zero vorticity 
in \cite{Hur2008solitary, GW2008solitary} and \cite{Wheeler2013solitary}. 
Specifically, for an arbitrary non-zero vorticity, one of the authors employed the generalized implicit function theorem of Nash-Moser type to construct a family of small-amplitude solitary water waves with super-critical wave speed near the KdV soliton. Moreover they are unique.

\subsection{The first-order approximation}\label{sec:solitary1}
We follow the same line of argument as in the previous section 
and develop the surface reconstruction procedure from the pressure at the fluid bed,
for small-amplitude solitary water wave with vorticity near the KdV soliton. 
We therefore assume that (see \cite{Hur2008solitary}, for instance)
\begin{equation}\label{E:solitary-parameter}
\delta = \sqrt{\epsilon}\ll 1
\end{equation}
and $u_j$'s, $v_j$'s, $p_j$'s and $\eta_j$'s, $j=1,2,3,\dots$, in \eqref{def:uv-scale}-\eqref{def:P-scale} 
decay to zero as $x\to\pm\infty$. For an arbitrary distribution of vorticity, incidentally, 
all solitary water waves decay to zero exponentially fast at infinity; see \cite{Hur-symmS}, for instance.

Under this assumption, \eqref{E:Euler-scaled} and \eqref{E:top-scaled}, \eqref{E:bottom-scaled} 
at the leading order become:
\begin{equation}\label{E:solitary-Euler1}
\left\{\begin{split}
&(U-c)u_{1x}+U'v_1= -p_{1x},\quad&& \\
&p_{1y}=0\qquad&&\text{in}\quad 0<y<h_0, \\
&u_{1x}+v_{1y}= 0, \\
\end{split}\right.
\end{equation}
and
\begin{align}
&v_1=(U-c)\eta_{1x} \quad\text{and}\quad g\eta_1=p_1 &&\text{at}\quad y=h_0, \label{E:solitary-top1}\\
&v_1= 0	&&\text{at}\quad y = 0.\label{E:solitary-bottom1}
\end{align}

Small-amplitude solitary water waves with vorticity constructed in \cite{Hur2008solitary}, for instance,
solve \eqref{E:steady}, and hence \eqref{E:Euler-scaled}-\eqref{E:bottom-scaled}.
It follows from the Lyapunov-Schmidt reduction in \cite{Hur2008solitary} that as $\delta^2=\epsilon\to0$, 
they are approximated by the KdV soliton with errors in $O(\epsilon^2)$ in the real analytic function space
and the wave speed is approximated by that determined upon solving \eqref{E:Burns} with errors in $O(\epsilon^2)$.
For an arbitrary distribution of vorticity in a H\"older space, incidentally,
one of the authors in \cite{Hur2012regularity} proved that 
a solitary water wave in the corresponding H\"older space is real analytic. 
Throughout the section, $c$ means the critical wave speed. 
Numerical computations in \cite{deconinck2012relating}, for instance, indicate that
the maximum wave height does not suffer much from assuming the critical wave speed.

Suppose that the dynamic pressure (see \eqref{def:P-scale}) at the bed is prescribed to the order of $\epsilon$. 
We write that
\begin{equation}\label{def:solitary-b1}
p_1(x,0)=b_1(x),
\end{equation}
and assume that $b_1$ is smooth and decays to zero as $x\to\pm\infty$. 
Since the second equation in \eqref{E:solitary-Euler1} implies that 
the dynamic pressure to the leading order does not vary with the depth, 
it follows from the second equation in \eqref{E:solitary-top1} 
the hydrostatic approximation (see \eqref{def:hydrostatic})
\begin{equation}\label{E:solitary-eta1}
g\eta_1(x)=b_1(x).
\end{equation}
This furnishes an explicit formula relating $\eta_1$ and $p_1(\cdot,0)$. 
It is independent of the underlying shear flow, and hence it is likely to be a poor approximation. 
In the following subsection, we shall compute higher-order correction terms to \eqref{E:solitary-eta1}, 
which do incorporate the effects of vorticity. 

For future usefulness, we infer from the first and the last equations in \eqref{E:solitary-Euler1} that
\begin{equation}\label{E:Rayleigh2}
\Big( \frac{v_1}{U-c}\Big)_y= \frac{b_1'}{(U-c)^2}.
\end{equation}
We solve it by quadrature to arrive at that $v_1(x,y)= b_1'(x)(U(y)-c)fr(y)$, where 
\begin{equation}\label{def:Fr}
fr(y)= \int_0^y \frac{dz}{(U(z)-c)^2}.
\end{equation}
As a matter of fact, one can solve \eqref{E:Rayleigh2} analytically for all $U$, 
unlike the Rayleigh equations in the previous section.
Note in passing that $fr(h_0)$ is the inverse square of the Froude number. 
We then determine $u_1$ upon integrating the last equation in \eqref{E:solitary-Euler1} 
and using that $u_1$ vanishes at infinity. To summarize,
\begin{equation}\label{E:solitary-soln1}
\left\{\begin{split}
&u_1(x,y)= -b_{1}(x)\left(U'(y)fr(y)+\frac{1}{U(y)-c}\right), \\
&v_1(x,y)= b_1'(x)(U(y)-c)fr(y),\\
&p_1(x,y)= b_1(x),
\end{split}\right.
\end{equation}
where $fr$ is defined in \eqref{def:Fr}.

\subsection{Higher-order approximations}\label{sec:solitary2}
We continue to assume \eqref{E:solitary-parameter}. To proceed, 
\eqref{E:Euler-scaled} and \eqref{E:top-scaled}, \eqref{E:bottom-scaled} at the order of $\epsilon$ become:
\begin{equation}\label{E:solitary-Euler2}\left\{
\begin{split}
&(U-c)u_{2x}+U'v_2+u_1u_{1x}+v_1u_{1y}= -p_{2x},\qquad && \\
&(U-c)v_{1x}= -p_{2y} &&\text{in}\quad 0<y<h_0, \\
&u_{2x}+v_{2y}= 0 &&
\end{split}\right.
\end{equation}
and 
\begin{align}
&v_2=(U-c)\eta_{2x}+u_1\eta_{1x} \quad\text{and}\quad g\eta_2=p_2	
\quad&&\text{at}\quad	y=h_0, \label{E:solitary-top2}\\
&v_2= 0	&&\text{at}\quad	y = 0. \label{E:solitary-bottom2}
\end{align}

Suppose that the dynamic pressure at the fluid bed is prescribed to the order of $\epsilon^2$. We write that
\begin{equation}\label{def:solitary-b2}
p_2(x,0)=b_2(x),
\end{equation}
and we continue to assume that $b_2$ is smooth and decays to zero as $x\to\pm\infty$. 
Integrating the second equation in \eqref{E:solitary-Euler2},
we use \eqref{E:solitary-soln1} and \eqref{def:solitary-b2} to arrive at that 
\[
p_2(x,y)=-b_1''(x)\int_0^y(U(z)-c)^2fr(z)~dz+b_2(x),
\]
where $fr$ is defined in \eqref{def:Fr}.
It then follows from the second equation in \eqref{E:solitary-top2} that
\begin{equation}\label{E:solitary-eta2}
g\eta_2(x)=-b_1''(x)\int_0^{h_0}(U(y)-c)^2fr(y)~dz+b_2(x).
\end{equation}
This furnishes an explicit formula relating $\eta_2$ and $p_1(\cdot,0)$ (see \eqref{def:solitary-b1}) and $p_2(\cdot,0)$, 
provided with the background shear flow and the wave speed satisfying \eqref{E:Burns}.
It incorporates the effects of the vorticity through the solution of the Rayleigh equation \eqref{E:Rayleigh2}.

We may repeat the argument in the previous subsection to find $u_2$ and $v_2$. 
Specifically, we infer from the first and the last equations in \eqref{E:solitary-Euler2} that 
\[
\Big( \frac{v_2}{U-c}\Big)_y= \frac{p_{2x}+u_1u_{1x}+v_1u_{1y}}{(U-c)^2},
\]
where $p_2$ is determined above and $u_1$ and $v_1$ are in \eqref{E:solitary-soln1}.
We determine $v_2$ upon integrating this. 
We then determine $u_2$ upon integrating the last equation in \eqref{E:solitary-Euler2} 
and using that $u_2$ vanishes at infinity. To summarize,
\begin{equation}\label{E:solitary-soln2}
\left\{\begin{split}
&u_{2}(x,y)= -U'(y)\int_0^y \int_{-\infty}^x \frac{p_{2x}+u_1u_{1x}+v_1u_{1y}}{(U-c)^2}(w,z)~dwdz\\
&\hspace*{100pt}-\int_{-\infty}^x\frac{p_{2x}+u_1u_{1x}+v_1u_{1y}}{U-c}(w,y)~dw,\\
&v_2(x,y)= (U(y)-c)\int_0^y \frac{p_{2x}+u_1u_{1x}+v_1u_{1y}}{(U-c)^2}(x,z)~dz, \\
&p_2(x,y)=-b_1''(x)\int_0^y(U(z)-c)^2fr(z)~dz+b_2(x),
\end{split}\right.
\end{equation}
where $fr$ is in \eqref{def:Fr} and $u_1$ and $u_2$ are in \eqref{E:solitary-soln1}.

Continuing, \eqref{E:Euler-scaled} and \eqref{E:top-scaled}, \eqref{E:bottom-scaled} 
at the order of $\epsilon^2$ become:
\begin{equation}\label{E:solitary-Euler3}\left\{
\begin{split}
&(U-c)u_{3x}+U'v_3+u_1u_{2x}+u_2u_{1x}+v_1u_{2y}+v_2u_{1y}	= -p_{3x},\\
&(U-c)v_{2x}+u_1v_{1x}+v_1v_{1y}= -p_{3y}, \\
&u_{3x}+v_{3y}= 0
\end{split}\right.
\end{equation}
in $0<y<h_0$ and
\begin{align}
&v_3=(U-c)\eta_{3x}+u_1\eta_{2x}+u_2\eta_{1x} \quad\text{and}\quad g\eta_3=p_3
\quad&&\text{at}\quad	y=h_0, \label{E:solitary-top3}\\
&v_3=0 &&\text{at}\quad y=0. \notag
\end{align}
Suppose that the dynamic pressure at the fluid bed is prescribed to the order of $\epsilon^3$:
\begin{equation}\label{def:solitary-b3}
p_3(x,0)=b_3(x),
\end{equation}
and we continue to assume that $b_3$ is smooth and decays to zero as $x\to\pm\infty$.
Integrating the second equation in \eqref{E:solitary-Euler3}, we use \eqref{def:solitary-b3} to arrive at that 
\[
 p_3(x,y)=-\int_0^y((U-c)v_{2x}+u_1v_{1x}+v_1v_{1y})(x,z)~dz+b_3(x),
\]
where $u_1$, $v_1$ and $v_2$ are in \eqref{E:solitary-soln1} and \eqref{E:solitary-soln2}.
Therefore, 
\begin{equation}\label{E:solitary-eta3}
g\eta_3(x)=-\int_0^{h_0}((U-c)v_{2x}+u_1v_{1x}+v_1v_{1y})(x,y)~dy+b_3(x).
\end{equation}
This furnishes an explicit formula relating $\eta_3$ and $p_j(\cdot,0)$, $j=1,2,3$, 
provided with the background shear flow and the wave speed satisfying \eqref{E:Burns}.
It incorporates the effects of vorticity through the solution of the Rayleigh equation \eqref{E:Rayleigh2}.

We may repeat the above argument and continue to higher-order approximations. 
As a matter of fact, one is able to derive explicit formulae relating $\eta_j$ as functions of $p_{j'}(\cdot,0)$
for each $j$ and for all $j'=1,2,\dots,j$. 
Expressions become quite complicated, however, and hence we do not pursue here.

\subsection{Examples}\label{sec:solitary.e.g.}

We illustrate the results in Section~\ref{sec:solitary1} and Section~\ref{sec:solitary2} 
by discussing some examples. Throughout the subsection, we assume for simplicity that
\[
b_2(x,0)=b_3(x,0)=0.
\]

\begin{example}[Zero vorticity]\label{e.g.:solitary0}\rm
In the case of $U\equiv0$, namely the zero vorticity, a straightforward calculation reveals that $fr(y)=y/c^2$, 
whence \eqref{E:solitary-soln1} and \eqref{E:solitary-eta1} becomes
\[
\left\{\begin{split}
&u_1(x,y)=\frac{b_1(x)}{c},\\
&v_1(x,y)=-\frac{b_1'(x)}{c}y,\\
&p_1(x,y)=b_1(x),
\end{split}\right.
\]
and $g\eta_1(x)=b_1(x)$. Of course, this represents the hydrostatic approximation \eqref{def:hydrostatic}. 
To proceed, \eqref{E:solitary-soln2} and \eqref{E:solitary-eta2} become
\[
\left\{\begin{split}
&u_2(x,y)=-\frac12\frac{b_1''(x)}{c}y^2+\frac12\frac{b_1^2(x)}{c^3}\\
&v_2(x,y)= \frac16\frac{b_1'''(x)}{c}y^3-\frac{b_1(x)b_1'(x)}{c^3}y\\
&p_2(x,y)=-\frac{1}{2}b_1''(x)y^2\\
\end{split}\right.\]
and
\begin{equation}\label{E:solitary0eta2}
g\eta_2(x)=-\frac{1}{2} h_0^2 b_1''(x).
\end{equation}
Continuing, \eqref{E:solitary-eta3} becomes
\begin{equation}\label{E:solitary0eta3}
g\eta_3(x)=\frac{1}{24}h_0^4b_1^{(4)}(x)-\frac{h_0b_1'(x)^2}{c}.
\end{equation}

We compare the results to asymptotic approximations of the exact formula in \cite{OVDH2012}, for instance.
The results in \cite[Section~4.2]{OVDH2012} may be written, in the dimensionless variables, abusing notation, as
\[
\left\{\begin{split}
&\eta_1(x)=b_1(x), \\
&\eta_2(x)=-\frac12b_1''(x), \\
&\eta_3(x)=\frac{1}{24}b_1^{(4)}(x)-b_1(x)b_1''(x)-\frac12b_1'(x)^2\Big(c^2+\frac{1}{c^2}\Big).
\end{split}\right.
\]
\end{example}
The second equation agrees with \eqref{E:solitary0eta2} up to physical constants 
and the last equation agrees with \eqref{E:solitary0eta3} up to physical constants except the middle term. 
Note that $c\approx 1$ in the non-dimensionalization.

\begin{example}[Constant vorticity]\rm
Let $U(y)=\gamma y$, $0\leq y\leq h_0$, for some constant $\gamma$.
(More generally, one may take $U(y)=\gamma y+U_0$ for some constant $U_0$, 
but $U_0$ may be absorbed into the wave speed.) This models the constant vorticity $-\gamma$. 

A straightforward calculation reveals that 
\[
fr(y)=\frac1c\frac{y}{c-\gamma y},
\]
whence \eqref{E:solitary-soln1} and \eqref{E:solitary-eta1} become
\[
\left\{\begin{split}
&u_1(x,y)=\frac{b_1(x)}{c},\\
&v_1(x,y)=-\frac{b_1'(x)}{c}y,\\
&p_1(x,y)=b_1(x),
\end{split}\right.
\]
and $g\eta_1=b_1(x)$. To proceed, \eqref{E:solitary-eta2} becomes
\[
g\eta_2(x)=\frac{b_1''(x)}{c}\Big(\frac13\gamma h_0^3-\frac12ch_0^2\Big).
\]
Other solutions may be computed explicitly using \eqref{E:solitary-soln2} and \eqref{E:solitary-eta3}.
Expressions are quite complicated, however, and we do not record them here.
\end{example}

\begin{example}[Poiseuille flows]\rm
Let $U(y)=h_0^2-y^2$, $0\leq y\leq h_0$. This models Poiseuille flows. 
Note that $c>\max_{0\leq y\leq h_0}U(y)=h_0^2$ to guarantee that solitary water waves exist.

One may evaluate \eqref{def:Fr} using Mathematica to find that
\[
fr(y)=\frac{1}{2\sqrt{c-h_0^2}^3}\Big(\frac{\sqrt{c-h_0^2}y}{y^2+c-h_0^2}+\tanh^{-1}\Big(\frac{y}{\sqrt{c-h_0^2}}\Big)\Big)
\]
and evaluate \eqref{E:solitary-eta2} to find that
\[\begin{split}
g\eta_2(x)=-\frac{b_1''(x)}{30(c-h_0^2)^{3/2}}\Big(&h_0(15 c^2-20 ch_0^2+8 h_0^4)\tan ^{-1}\Big(\frac{h_0}{\sqrt{c-h_0^2}}\Big) \\&+\sqrt{c-h_0^2}\Big(h_0^2(4 c-h_0^2)+4(c-h_0^2)^2 \log\Big(\frac{c-h_0^2}{c}\Big)\Big)\Big).
\end{split}\]
Other solutions may be computed explicitly. Expressions are quite complicated and we do not record them here.
\end{example}

\begin{appendix}

\section{The zero vorticity case, revisited}\label{sec:OVDH2012}

We discuss the exact formula in \cite{OVDH2012} relating the Stokes wave with zero vorticity and the pressure at the bed.

Under the assumption that the flow is irrotational, one may rewrite the governing equations of the water wave problem
in terms of the velocity potential $\phi(x,y;t)$ and the surface elevation $\eta(x;t)$ from the undisturbed depth: 
\begin{equation}\label{ww-problem}
\left\{\begin{split}
&\Delta \phi = 0 &&\text{in}\quad 0<y<h_0+\eta(x;t),\\
&\phi_t + \frac{1}{2} (\phi_x^2 + \phi_y^2) =p &&\text{in}\quad0<y<h_0+\eta(x;t), \\
&\phi_y = \eta_t + \phi_x \eta_x\quad\text{and}\quad p = g\eta\qquad &&\text{at}\quad y = h_0+\eta(x;t),\\
&\phi_y = 0 &&\text{at}\quad y=0.
\end{split}\right. 
\end{equation}
Recall that $p$ is the dynamic pressure, measuring the departure from the hydrostatic pressure.

Let $q(x;t) = \phi(x,h_0+\eta(x;t);t)$ represent the trace of the velocity potential at the fluid surface $y=\eta(x;t)$. 
We appeal to the chain rule and use the former of the third equations in \eqref{ww-problem} to show that
\[
\phi_x = \frac{q_x-\eta_x \eta_t}{1+\eta_x^2} \quad\text{and}\quad \phi_y = \frac{\eta_t + \eta_x q_x}{1+\eta_x^2} 
\qquad\text{at}\quad y=h_0+\eta(x;t).
\]
Similarly,
\[
\phi_t = q_t-\frac{\eta_t(\eta_t+\eta_x q_x)}{1+\eta_x^2}\qquad\text{at}\quad y=h_0+\eta(x;t).
\]
In the moving coordinate frame, where $q$ and $\eta$ are functions of $x-ct$,
\[
\phi_x = \frac{q_x+c \eta_x^2}{1+\eta_x^2},\qquad \phi_y = \frac{(q_x - c)\eta_x}{1+\eta_x^2}
\quad\text{and}\quad \phi_t = -c q_x +\frac{c \eta_x^2(q_x-c)}{1+\eta_x^2}
\]
at $y=h_0 +\eta(x)$. After substitution,
the second equation and the latter of the third equations in \eqref{ww-problem} imply that
\[
q_x^2-2cq_x-c^2 \eta_x^2 + 2g\eta(1+\eta_x^2)=0.
\]
Therefore $q_x = c \pm \sqrt{(c-2g\eta)(1+\eta_x^2)}$. 
We choose the $-$ sign to guarantee that $u-c<0$ throughout the fluid region. Consequently
\begin{equation}\label{E:exact-eta}
\phi_x(x,h_0+\eta(x))=c-\sqrt{\frac{c^2-2g\eta(x)}{1+\eta_x^2(x)}}.
\end{equation}
Restricting the second equation in \eqref{ww-problem} at $y=0$, moreover, 
we arrive in the moving coordinates frame at that
\[
-c\phi_x+\frac{1}{2} \phi_x^2=-p \qquad\text{at}\quad y=0.
\]
Therefore $\phi_x(x,0)=c\pm\sqrt{c^2-2p(x,0)}$. 
We choose the $-$ sign to guarantee that $u-c<0$ throughout the fluid region.
To recapitulate, in the surface reconstruction from the pressure at the bed for Stokes waves with zero vorticity, 
one solves the boundary value problem for the Laplace equation:
\begin{equation}\label{E:periodic0prob}
\left\{\begin{split}
&\Delta \phi = 0\qquad \text{in}\quad y>0,\\
&\phi_x = c-\sqrt{c^2-2p(x,0)}\quad\text{and}\quad \phi_y = 0\qquad \text{at}\quad y=0;
\end{split}\right.
\end{equation}
one then uses \eqref{E:exact-eta} to determine $\eta$ in terms of $p(x,0)$.

Let's write that
\[
\phi(x,y)=\sum_{-\infty}^{\infty} \phi_n(y) e^{inx}\quad\text{and}\quad
c-\sqrt{c^2-2p(x,0)}=\sum_{-\infty}^{\infty} q_ne^{inx}.
\]
In other words, 
\[
q_n=\frac{1}{2\pi}\int^{2\pi}_0 (c-\sqrt{c^2-2p(x,0)})e^{inx}~dx.
\]
For each $n\neq 0$, note that \eqref{E:periodic0prob} leads to 
the Cauchy problem for the linear second-order constant-coefficient ODE:
\[
\left\{\begin{split}
&\phi_n''-n^2\phi_n=0\qquad \text{for}\quad y>0,\\
&\phi_n'(0)=0\quad\text{and}\quad in\phi_n(0)=q_n,
\end{split}\right. 
\]
whence
\[
\phi(x,y)=\sum_{n\neq 0} -iq_n\cosh(ny) \frac{e^{inx}}{n}
\]
up to addition by a constant. Note that $q_0=0$ and $\phi_0$ is an arbitrary constant.
It then follows from \eqref{E:exact-eta} that
\begin{equation}\label{E:exact-periodic}
\begin{split}
c - \sqrt{\frac{c^2 - 2g\eta}{1+\eta_x^2}}=&\sum_{n\neq0} q_n \cosh(n(h_0+\eta(x))) e^{inx}\\
=&\sum_{n\neq0}\Big(\int^{2\pi}_0 (c-\sqrt{c^2-2p(x,0)})e^{inx}~dx\Big)\cosh(n(h_0+\eta(x)))e^{inx}.
\end{split}
\end{equation}
This furnishes an implicit formula relating $\eta$ and $p(\cdot,0)$, provided with a suitable wave speed.
It agrees with (22) in \cite{OVDH2012}, for which the sum on the right side ranges over all integers, since $q_0=0$.
Nevertheless, we emphasize that the zeroth Fourier mode is excluded, 
which is important below in the derivation of asymptotic approximations.

To proceed, we assume that $\epsilon\ll 1$ and that
\[
\eta= \epsilon \eta_1+ \epsilon^2 \eta_2+ \cdots\quad\text{and}\quad 
p(\cdot,0)=\epsilon p_1.
\]
Substituting these, we expand the left side of \eqref{E:exact-periodic} in the Taylor fashion to write that 
\begin{align}\label{E:LHS}
c-\sqrt{\frac{c^2-2g\eta}{1+\eta_x^2}}
=&c-\sqrt{(c^2-2g\epsilon\eta_1-2g\epsilon^2 \eta_2-\cdots)(1-\epsilon^2 \eta_{1x}^2+\cdots)}\notag\\
=&c-\sqrt{c^2-2g\epsilon\eta_1-2g\epsilon^2\eta_2-\epsilon^2 c^2 \eta_{1x}^2+\cdots}\notag\\
=&c-c\sqrt{1-2\Big(\frac{g}{c^2}\epsilon\eta_1
+\frac{g}{c^2}\epsilon^2\eta_2+\epsilon^2 \eta_{1x}^2+\cdots\Big)}\notag\\
=&c-c\Big(1-\Big(\frac{g}{c^2}\epsilon\eta_1+\frac{g}{c^2}\epsilon^2\eta_2+\epsilon^2 \eta_{1x}^2\Big)
-\frac{1}{2}\frac{g^2}{c^4}\epsilon^2 \eta_1^2+\cdots\Big)\notag\\
=&\frac{g}{c}\epsilon\eta_1+\frac{g}{c}\epsilon^2\eta_2+\frac{1}{2}c\epsilon^2 \eta_{1x}^2+\frac{g^2}{2 c^3}\epsilon^2\eta_1^2+\cdots.
\end{align}
Similarly,
\begin{align}
\cosh(n(h_0+\eta))=&\cosh(n h_0+n\epsilon \eta_1+n\epsilon^2 \eta_2+\cdots)\notag\\
=& \cosh(n h_0) +\sinh(n h_0)(n\epsilon \eta_1+n\epsilon^2 \eta_2+...)+\frac{1}{2}\cosh(n^2\epsilon^2 \eta_1^2+\cdots).\label{E:RHS1}
\end{align}
and
\begin{align}
c-\sqrt{c^2-2b}=&c-\sqrt{c^2-2\epsilon p_1+\cdots} \notag\\
=&c-c\sqrt{1-2\frac{1}{c^2}\epsilon p_1+\cdots}
=\frac{1}{c}\epsilon p_1+\frac{1}{2 c^3}\epsilon^2 p_1^2+\cdots.\label{E:RHS2}
\end{align}

Let's write that
\[
p_1(x)=\epsilon \sum_{n=-\infty}^\infty b_ne^{inx}.
\] 
Substituting \eqref{E:LHS}-\eqref{E:RHS2} into \eqref{E:exact-periodic}, at the order of $\epsilon$, we gather that
\[
g \eta_1(x)= \sum_{n\neq0}\cosh(n h_0) b_n e^{inx}.
\]
This agrees with (51) in \cite{OVDH2012} if it is written in the dimensional variables. 
Continuing, at the order of $\epsilon^2$, we gather that
\begin{equation}\label{E:OVDH2}
\frac{g}{c} \eta_2 + \frac{c}{2} \eta_{1x}^2 + \frac{g^2}{c^3} \eta_1^2
=\sum_{n\neq0}\cosh(nh_0)\frac12\frac{p_1^2}{c^3}+\eta_1\sum_{n\neq0}n\sinh(nh_0)\frac{p_1}{c}.
\end{equation}
This is similar to (52) in \cite{OVDH2012} in the dimensional variables, 
but the convolution sums must ranges over nonzero integers in the Stokes wave setting. 
To illustrate and to compare the results with those in Example~\ref{e.g.:periodic0}, we furthermore assume that 
\[
p_1(x)=b\cos(x).
\] 
We then find, instead, that
\begin{equation}\label{E:OVDH1'}
g\eta_1(x)=\cosh(h_0)b\cos(x).
\end{equation} 
To proceed, \eqref{E:OVDH2} becomes
\begin{align*}
g\eta_2=&-\frac12c^2\frac{b^2}{g^2}\cosh^2(h_0)\sin^2(x)-\frac12\frac{g^2}{c^2}\frac{b^2}{g^2}\cosh^2(h_0)\cos^2(x)\\
&+\sum_{n\neq0}\cosh(nh_0)\frac12\frac{b^2}{c^2}\cos^2(x)
+\frac{b}{g}\cosh(h_0)\cos(x)\sum_{n\neq0}n\sinh(nh_0)b\cos(x) \\
=&-\frac12\frac{b^2}{g}\sinh(h_0)\cosh(h_0)\frac{1-\cos(2x)}{2}
-\frac12\frac{b^2}{g}\frac{\cosh^3(h_0)}{\sinh(h_0)}\frac{1+\cos(2x)}{2} \\
&+\frac12\frac{b^2}{g}\frac{\cosh(h_0)}{\sinh(h_0)}\cosh(2h_0)\frac12\cos(2x)
+\frac{b^2}{g}\cosh(h_0)\sinh(h_0)\frac{1+\cos(2x)}{2}\\
=& \left(\frac{b}{2c} \right)^2 (-1+4\sinh^2(h_0) \cos(2x)).
\end{align*}
The second and third equalities use \eqref{E:c0}, where $k=1$. 
The constant term on the right side becomes zero if the sum ranges over all integers. 

\end{appendix}

\subsection*{Acknowledgements}
VMH is supported by the National Science Foundation grant CAREER DMS-1352597, an Alfred P. Sloan research fellowship, and an Arnold O. Beckman research award RB14100, the Center for Advanced Study Beckman fellowship at the University of Illinois at Urbana-Champaign. 
MRL is supported through an Arnold O. Beckman research award RB14100 at the University of Illinois at Urbana-Champaign.

\bibliographystyle{amsalpha}
\bibliography{transferBib}

\end{document}